\documentclass{article}
\usepackage{amsmath}
\usepackage{amsfonts}
\usepackage[ansinew]{inputenc}

\usepackage{harvard}

\newcommand{\EE}{\mathbb{E}}
\newcommand{\PP}{\mathbb{P}}
\newcommand{\RR}{\mathbb{R}}
\newcommand{\VV}{\mathbb{V}}

\newcommand{\cov}{\operatorname{cov}}
\newcommand{\corr}{\operatorname{corr}}

\newcommand{\He}{\operatorname{He}}

\newtheorem{theorem}{Theorem}[section]

\numberwithin{equation}{section}

\begin{document}

\title{The Bivariate Normal Copula}
\author{Christian Meyer\footnote{DZ BANK AG, Platz der Republik, D-60265 Frankfurt. The opinions or recommendations expressed in this article are those of the author and are not representative of DZ BANK AG.} \footnote{E-Mail: {\tt Christian.Meyer@dzbank.de}}}
\date{\today}
\maketitle

\begin{abstract}
We collect well known and less known facts about the bivariate normal distribution and translate them into copula language. In addition, we prove a very general formula for the bivariate normal copula, we compute Gini's gamma, and we provide improved bounds and approximations on the diagonal.
\end{abstract}

\section{Introduction}

When it comes to modelling dependent random variables, in practice one often resorts to the multivariate normal distribution, mainly because it is easy to parameterize and to deal with. In recent years, in particular in the quantitative finance community, we have also witnessed a trend to separate marginal distributions and dependence structure, using the copula concept. The normal (or Gauss, or Gaussian) copula has even come to the attention of the general public due to its use in the valuation of structured products and the decline of these products during the financial crisis of 2007 and 2008.

The multivariate normal distribution has been studied since the 19th century. Many important results have been published in the 1950s and 1960s. Nowadays, quite some effort is spent on rediscovering some of them.

In this note we will proceed as follows. We will concentrate on the bivariate normal distribution because it is the most important case in practice (applications in quantitative finance include pricing of options and estimation of asset correlations; the impressive list by \citeasnoun{BL} mentions applications in agriculture, biology, engineering, economics and finance, the environment, genetics, medicine, psychology, quality control, reliability and survival analysis, sociology, physical sciences and technology). We will give an extensive view on its properties. Everything will be reformulated in terms of the associated copula.

We will also provide new results (at least if the author has not been rediscovering himself...), including a very general formula implying other well known formulas, a derivation of Gini's gamma, and improved and very simple bounds for the diagonal.

For collections of facts on the bivariate (or multivariate) distribution we refer to the books of \citeasnoun{BL}, of \citeasnoun{KBJ}, and of \citeasnoun{PR}, and to the survey article of \citeasnoun{Gupta} with its extensive bibliography \cite{GuptaBib}. For theory on copulas we refer to the book of \citeasnoun{Nelsen}.

We will use the symbols $\PP$, $\EE$ and $\VV$ for probabilities, expectation values and variances.


\section{Definition and basic properties}
\label{sec_basic}

Denote by
\[
\varphi(x) := \frac{1}{\sqrt{2\pi}} \exp\left(-\frac{x^2}{2}\right),
\qquad
\Phi(h) := \int_{-\infty}^h \varphi(x)\;dx
\]
the density and distribution function of the standard normal distribution, and by
\begin{align*}
\varphi_2(x,y;\varrho) & := \frac{1}{2\pi\sqrt{1-\varrho^2}}
\exp\left(-\frac{x^2-2\varrho xy + y^2}{2(1-\varrho^2)}\right),\\
\Phi_2(h,k;\varrho) & := \int_{-\infty}^{h} \int_{-\infty}^{k}
\varphi_2(x,y;\varrho) \; dy \; dx,
\end{align*}
the density and distribution function of the bivariate standard normal distribution with correlation parameter $\varrho\in(-1,1)$. The bivariate normal (or Gauss, or Gaussian) copula with parameter $\varrho$ is then defined by application of Sklar's theorem, cf. Section 2.3 of \cite{Nelsen}:
\begin{equation}
C(u,v;\varrho) := \Phi_2\left(\Phi^{-1}(u),\Phi^{-1}(v);\varrho\right)
\label{eq_def}
\end{equation}
For $\varrho\in\{-1,1\}$, the correlation matrix of the bivariate standard normal distribution becomes singular. Nevertheless, the distribution, and hence the normal copula, can be extended continuously. We may define
\begin{align}
C(u,v,-1) & := \lim_{\varrho\longrightarrow -1} C(u,v;\varrho) = \max(u+v-1,0),\label{eq_min}\\
C(u,v,+1) & := \lim_{\varrho\longrightarrow +1} C(u,v;\varrho) = \min(u,v).\label{eq_max}
\end{align}
Hence $C(\cdot,\cdot;\varrho)$, for $\varrho\longrightarrow -1$, approaches the lower Fr\'echet bound,
\[
W(u,v) := \max(u+v-1,0),
\]
and, for $\varrho\longrightarrow 1$, approaches the upper Fr\'echet bound,
\[
M(u,v) := \min(u,v).
\]
Furthermore, we have
\begin{equation}
C(u,v;0) = u\cdot v =:\Pi(u,v),\label{eq_zero}
\end{equation}
the independence copula. What happens in between can be described using the following differential equation derived by \citeasnoun{Plackett}:
\begin{equation}
\frac{d}{d\varrho} \Phi_2(x,y;\varrho) = \varphi_2(x,y;\varrho)
= \frac{d^2}{dx \; dy} \Phi_2(x,y;\varrho)
\end{equation}
We find
\begin{equation}
C(u,v;\varrho) - C(u,v;\sigma) = I(u,v;\sigma,\varrho)
:= \int_{\sigma}^{\varrho} \varphi_2\left(\Phi^{-1}(u),\Phi^{-1}(v);r\right) \; dr,
\end{equation}
and in particular,
\begin{align}
C(u,v;\varrho)
& = W(u,v) + I(u,v;-1,\varrho),\label{eq_plackett_min}\\
& = \Pi(u,v) + I(u,v;0,\varrho),\label{eq_plackett_ind}\\
& = M(u,v) - I(u,v;\varrho,1).\label{eq_plackett_max}
\end{align}
In other words, the bivariate normal copula allows comprehensive total concordance ordering with respect to $\varrho$, cf. Section 2.8 of \cite{Nelsen}:
\[
W(\cdot,\cdot) = C(\cdot,\cdot;-1)
\prec C(\cdot,\cdot;\varrho)
\prec C(\cdot,\cdot;\sigma)
\prec C(\cdot,\cdot;1) = M(\cdot,\cdot)
\]
for $-1\leq\varrho\leq\sigma\leq 1$, i.e., for all $u$, $v$,
\[
W(u,v) = C(u,v;-1) \leq
C(u,v;\varrho) \leq C(u,v;\sigma) \leq C(u,v,1) = M(u,v).
\]
By plugging $u=v=\frac{1}{2}$ into (\ref{eq_plackett_ind}) we obtain
\begin{equation}
C\left(\textstyle\frac{1}{2},\frac{1}{2};\varrho\right)
= \frac{1}{4} + \int_0^{\varrho} \frac{1}{2\pi\sqrt{1-r^2}}\;dr
= \frac{1}{4} + \frac{1}{2\pi}\arcsin(\varrho),
\label{eq_center}
\end{equation}
a result already known to \citeasnoun{Stieltjes}.

\citeasnoun{Mehler} and \citeasnoun{Pearson1}, among other authors, obtained the so-called tetrachoric expansion in $\varrho$:
\begin{equation}
C(u,v;\varrho) = uv + \varrho\left(\Phi^{-1}(u)\right) \varrho\left(\Phi^{-1}(v)\right) \sum_{k=0}^{\infty} \He_k\left(\Phi^{-1}(u)\right) \He_k\left(\Phi^{-1}(v)\right) \frac{\varrho^{k+1}}{(k+1)!}
\label{eq_tetra}
\end{equation}
where
\[
\He_k(x) = \sum_{i=0}^{[k/2]} \frac{k!}{i!(k-2i)!}\left(-\frac{1}{2}\right)^i x^{k-2i}
\]
are the Hermite polynomials.

The bivariate normal copula inherits the symmetries of the bivariate normal distribution:
\begin{align}
C(u,v;\varrho) & = C(v,u;\varrho)\label{eq_swap}\\
& = u - C(u,1-v;-\varrho)\label{eq_reflect_v}\\
& = v - C(1-u,v;-\varrho)\label{eq_reflect_u}\\
& = u + v - 1 + C(1-u,1-v;\varrho)\label{eq_reflect_uv}
\end{align}
Here, (\ref{eq_swap}) is a consequence of exchangeability, and (\ref{eq_reflect_uv}) is a consequence of radial symmetry, cf. Section 2.7 of \cite{Nelsen}.

In the following sections we will discuss formulas for the bivariate normal copula, its numerical evaluation, bounds and approximations, measures of concordance, and univariate distributions related to the bivariate normal copula. It will be convenient to assume $\varrho\not\in\{-1,0,1\}$ unless explicitly stated otherwise (cf. (\ref{eq_min}), (\ref{eq_zero}), (\ref{eq_max}) for the simple formulation in the missing cases).


\section{Formulas}
\label{sec_formulas}

If $(X,Y)$ are bivariate standard normally distributed with correlation $\varrho$ then $Y$ conditional on $X=x$ is normally distributed with expectation $\varrho x$ and variance $1-\varrho^2$. This translates into the following formulas for the conditional distributions of the bivariate normal copula:
\begin{align}
\frac{\partial}{\partial u} C(u,v;\varrho)
& = \Phi\left(\frac{\Phi^{-1}(v)-\varrho\cdot\Phi^{-1}(u)}{\sqrt{1-\varrho^2}}\right),\label{eq_cond1}\\
\frac{\partial}{\partial v} C(u,v;\varrho)
& = \Phi\left(\frac{\Phi^{-1}(u)-\varrho\cdot\Phi^{-1}(v)}{\sqrt{1-\varrho^2}}\right)\label{eq_cond2}
\end{align}
The copula density is given by:
\begin{align}
c(u,v;\varrho) & := \frac{\partial^2}{\partial u \;\partial v} C(u,v;\varrho)
= \frac{\varphi_2\left(\Phi^{-1}(u),\Phi^{-1}(v);\varrho\right)}{\varphi\left(\Phi^{-1}(u)\right)\varphi\left(\Phi^{-1}(v)\right)}\\
& \phantom{:}= \frac{1}{\sqrt{1-\varrho^2}}
\exp\left(\frac{2\varrho\Phi^{-1}(u)\Phi^{-1}(v)-\varrho^2\left(\Phi^{-1}(u)^2+\Phi^{-1}(v)^2\right)}{2(1-\varrho^2)}\right)\nonumber
\end{align}
In terms of its copula density, the bivariate normal copula can be written as
\begin{equation}
C(u,v;\varrho) = \int_0^u \int_0^v c(s,t;\varrho)\;dt\;ds.
\label{eq_11c}
\end{equation}

Now let $\alpha,\beta,\gamma\in(-1,1)$ with $\alpha\beta\gamma =\varrho$. Then the following holds:
\begin{align}
& C(u,v;\varrho)\nonumber\\
& = \int_{-\infty}^{\infty} \int_{-\infty}^{\infty} \Phi\left(\frac{\Phi^{-1}(u)-\alpha\cdot x}{\sqrt{1-\alpha^2}}\right)
\Phi\left(\frac{\Phi^{-1}(v)-\beta\cdot y}{\sqrt{1-\beta^2}}\right)
\varphi_2(x,y;\gamma) \;dy \;dx\label{eq_abc}\\
& = \int_0^1 \int_0^1 \Phi\left(\frac{\Phi^{-1}(u)-\alpha\cdot \Phi^{-1}(s)}{\sqrt{1-\alpha^2}}\right)
\Phi\left(\frac{\Phi^{-1}(v)-\beta\cdot \Phi^{-1}(t)}{\sqrt{1-\beta^2}}\right) c(s,t;\gamma) \;dt \;ds\nonumber\\
& = \int_0^1 \int_0^1
\left( \frac{\partial}{\partial s} C(u,s;\alpha)\right)
\left( \frac{\partial}{\partial t} C(t,v;\beta)\right)
c(s,t;\gamma) \;dt \;ds\nonumber
\end{align}
The right hand side of (\ref{eq_abc}) has occurred in credit risk modelling but without the link to the bivariate normal copula, cf. \cite{BO}. A proof in that context will be provided in Section \ref{sec_model}.

Other formulas for $C(u,v;\varrho)$ are obtained by carefully studying the limit as some of the variables $\alpha$, $\beta$, $\gamma$ are approaching the value one. The interesting cases, not regarding symmetry, are listed in Table \ref{table_spec}.

\begin{table}[htbp]
\begin{center}
\begin{tabular}{cccc}
\hline
$\alpha$ & $\beta$ & $\gamma$ & reference\\
\hline
$=\alpha$ & $=\beta$ & $=\gamma$ & (\ref{eq_abc})\\
$=1$ & $=\beta$ & $=\gamma$ & (\ref{eq_1bc})\\
$=\alpha$ & $=\beta$ & $=1$ & (\ref{eq_ab1})\\
$=1$ & $=\varrho$ & $=1$ & (\ref{eq_1b1})\\
$=1$ & $=1$ & $=\varrho$ & (\ref{eq_11c})\\
\hline
\end{tabular}
\end{center}
\caption{Limiting cases of (\ref{eq_abc})}
\label{table_spec}
\end{table}

By approaching $\alpha=1$ in (\ref{eq_abc}) we obtain:
\begin{align}
C(u,v;\varrho)
& = \int_{-\infty}^{\Phi^{-1}(u)} \int_{-\infty}^{\infty}
\Phi\left(\frac{\Phi^{-1}(v)-\beta\cdot y}{\sqrt{1-\beta^2}}\right)
\varphi_2(x,y;\gamma) \;dy \;dx\label{eq_1bc}\\
& = \int_0^u \int_0^1
\Phi\left(\frac{\Phi^{-1}(v)-\beta\cdot \Phi^{-1}(t)}{\sqrt{1-\beta^2}}\right) c(s,t;\gamma) \;dt \;ds\nonumber\\
& = \int_0^u \int_0^1 \left( \frac{\partial}{\partial t} C(u,t;\beta)\right)
c(s,t;\gamma) \;dt \;ds\nonumber
\end{align}
However, Equation (\ref{eq_1bc}) may be considered rather unattractive. By approaching $\gamma=1$ in (\ref{eq_abc}) instead we obtain a more interesting formula:
\begin{align}
C(u,v;\varrho)
& = \int_{-\infty}^{\infty}
\Phi\left(\frac{\Phi^{-1}(u)-\alpha\cdot z}{\sqrt{1-\alpha^2}}\right)
\Phi\left(\frac{\Phi^{-1}(v)-\beta\cdot z}{\sqrt{1-\beta^2}}\right)
\varphi(z)\;dz\label{eq_ab1}\\
& = \int_0^1
\Phi\left(\frac{\Phi^{-1}(u)-\alpha\cdot\Phi^{-1}(t)}{\sqrt{1-\alpha^2}}\right)
\Phi\left(\frac{\Phi^{-1}(v)-\beta\cdot\Phi^{-1}(t)}{\sqrt{1-\beta^2}}\right)\;dt\nonumber\\
& = \int_0^1
\left( \frac{\partial}{\partial t} C(t,v;\alpha)\right)
\left( \frac{\partial}{\partial t} C(u,t;\beta)\right)\;dt\nonumber
\end{align}
Equation (\ref{eq_ab1}) seems to have been discovered repeatedly, sometimes in more general context (e.g., multivariate, cf. \cite{SteckOwen}), sometimes for special cases. \citeasnoun{Gupta} gives credit to \citeasnoun{DS}.

By approaching $\alpha=\gamma=1$ in (\ref{eq_abc}) we obtain:
\begin{align}
C(u,v;\varrho) & = \int_0^u
\Phi\left(\frac{\Phi^{-1}(v)-\varrho\cdot\Phi^{-1}(t)}{\sqrt{1-\varrho^2}}\right)\;dt
=\int_0^u \frac{\partial}{\partial t} C(t,v;\varrho)\;dt\label{eq_1b1}\\
& = \int_0^v
\Phi\left(\frac{\Phi^{-1}(u)-\varrho\cdot\Phi^{-1}(t)}{\sqrt{1-\varrho^2}}\right)\;dt
= \int_0^v \frac{\partial}{\partial t} C(u,t;\varrho)\;dt\label{eq_1b1v}
\end{align}
These formulas can also be derived from (\ref{eq_cond1}) and (\ref{eq_cond2}). Finally, by approaching $\alpha=\beta=1$ in (\ref{eq_abc}) we rediscover (\ref{eq_11c}).

In order to evaluate the bivariate normal distribution function numerically, \citeasnoun{Owen} defined the following very useful function, to which we will refer as Owen's $T$-function:
\begin{equation}
T(h,a) = \frac{1}{2\pi}\int_0^a \frac{\exp\left(-\frac{1}{2} h^2 (1+x^2)\right)}{1+x^2}\; dx
= \varphi(h) \int_0^a \frac{\varphi(hx)}{1+x^2}\; dx
\label{eq_T1}
\end{equation}
He proved that
\begin{equation}
C(u,v;\varrho) = \frac{u+v}{2} - T\left(\Phi^{-1}(u),\alpha_u\right) - T\left(\Phi^{-1}(v),\alpha_v\right) - \delta(u,v)
\label{eq_T2}
\end{equation}
where
\begin{equation}
\delta(u,v) :=
\begin{cases}
\frac{1}{2}, & \text{if}\quad u<\frac{1}{2},v\geq \frac{1}{2}\quad\text{or}\quad u\geq \frac{1}{2},v<\frac{1}{2}\\
0, & \text{else}
\end{cases}
\end{equation}
and
\begin{equation}
\alpha_u = \frac{1}{\sqrt{1-\varrho^2}}\left(\frac{\Phi^{-1}(v)}{\Phi^{-1}(u)}-\varrho\right),
\qquad
\alpha_v = \frac{1}{\sqrt{1-\varrho^2}}\left(\frac{\Phi^{-1}(u)}{\Phi^{-1}(v)}-\varrho\right).
\end{equation}
In particular, on the lines defined by $v=\frac{1}{2}$ and by $u=v$, the following formulas hold:
\begin{align}
C\left(u,\textstyle\frac{1}{2};\varrho\right) & = \frac{u}{2} - T\left(\Phi^{-1}(u),-\frac{\varrho}{\sqrt{1-\varrho^2}}\right),\label{eq_T_line}\\
C(u,u;\varrho) & = u - 2\cdot T\left(\Phi^{-1}(u),\sqrt{\frac{1-\varrho}{1+\varrho}}\right)
\end{align}
From (\ref{eq_T2}) and (\ref{eq_T_line}) we can derive the useful formula 
\begin{equation}
C(u,v;\varrho) = C\left(u,\textstyle\frac{1}{2};\varrho_u\right) + C\left(v,\textstyle\frac{1}{2};\varrho_v\right) - \delta(u,v),
\label{eq_uv2line}
\end{equation}
where
\begin{align*}
\varrho_u & = -\frac{\alpha_u}{\sqrt{1+\alpha_u^2}} = \sin(\arctan(-\alpha_u)),\\
\varrho_v & = -\frac{\alpha_v}{\sqrt{1+\alpha_v^2}} = \sin(\arctan(-\alpha_v)).
\end{align*}
On the diagonal $u=v$, (\ref{eq_uv2line}) reads:
\begin{equation}
C(u,u;\varrho) = 2\cdot C\left(u,\textstyle\frac{1}{2};\displaystyle-\sqrt{\frac{1-\varrho}{2}}\right)
\label{eq_diag2line}
\end{equation}
Inversion of (\ref{eq_diag2line}) using (\ref{eq_reflect_v}) gives:
\begin{equation}
C\left(u,\textstyle\frac{1}{2};\varrho\right) =
\begin{cases}
\frac{1}{2} C(u,u;1-2\varrho^2), & \varrho < 0,\\
u - \frac{1}{2}C(u,u;1-2\varrho^2), & \varrho > 0.
\end{cases}
\label{eq_line2diag}
\end{equation}

Applying (\ref{eq_1b1}) to (\ref{eq_diag2line}) we obtain, cf. also \cite{SteckOwen},
\begin{equation}
C(u,u;\varrho) = 2\cdot \int_{0}^{u} g(t;\varrho)\;dt
\label{eq_C_diag_int} 
\end{equation}
with
\begin{equation}
g(u;\varrho) := \Phi\left(\sqrt{\frac{1-\varrho}{1+\varrho}}\cdot\Phi^{-1}(u)\right).
\label{eq_g}
\end{equation}
We find
\begin{equation}
\frac{d}{du} C(u,u;\varrho) = 2\cdot g(u;\varrho).
\end{equation}
The function $g$ will become important in Section \ref{sec_bounds}. Note that if $U$, $V$ are uniformly distributed on $[0,1]$ with copula $C(\cdot,\cdot;\varrho)$ then
\[
g(u;\varrho) = \PP(V\leq u\,|\,U=u).
\]
Below we list some properties of $g$:
\begin{align}
\lim_{u\longrightarrow 0^+} g(u;\varrho) & = 0,\label{eq_g_lim_0}\\
\lim_{u\longrightarrow 1^-} g(u;\varrho) & = 1,\label{eq_g_lim_1}\\
g(\textstyle\frac{1}{2};\varrho) & = \textstyle\frac{1}{2},\\
g(1-u;\varrho) & = 1 - g(u;\varrho),\\
g(g(u;\varrho);-\varrho) & = u\label{eq_g_inv}
\end{align}
In particular, (\ref{eq_g_lim_0}) and (\ref{eq_g_lim_1}) show that the bivariate normal copula does not exhibit (lower or upper) tail dependence (cf. Section 5.2.3 of \citeasnoun{EFM}).

Substitution of $t=g(s;\varrho)$ in (\ref{eq_C_diag_int}) and application of (\ref{eq_g_inv}) lead to the identity, cf. also \cite{SteckOwen}:
\begin{equation}
C(u,u;\varrho) = 2u\cdot g(u;\varrho) - C\left(g(u;\varrho),g(u;\varrho);-\varrho\right)
\label{eq_gtrafo}
\end{equation}


\section{Numerical evaluation}
\label{sec_numerical}

The bivariate normal copula has to be evaluated numerically. To my knowledge, in the literature there is no direct approach. Hence for the time being we have to rely on the definition of the bivariate normal copula (\ref{eq_def}) and numerical evaluation of $\Phi_2$ and $\Phi^{-1}$. There are excellent algorithms available for evaluation of $\Phi^{-1}$, cf. \cite{Acklam}; the main problem is evaluation of the bivariate normal distribution function $\Phi_2$.

In the literature on evaluation of $\Phi_2$ there are basically two approaches: application of a multivariate method to the bivariate case, and explicit consideration of the bivariate case. For background on multivariate methods we refer to the recent book by \citeasnoun{BretzGenz}. In most cases, bivariate methods will be able to obtain the desired accuracy in less time. In the following we will provide an overview on the literature. We will concentrate on methods and omit references dealing with implementation only. Comparisons of different approaches in terms of accuracy and running time have been provided by numerous authors, e.g., \citeasnoun{AC}, \citeasnoun{TW}, and \citeasnoun{WangKen}.

Before the advent of widely available computer power, extensive tables of the bivariate normal distribution function had to be created. Using (\ref{eq_uv2line}) or similar approaches, the three-dimensional problem (two variables and the correlation parameter) was reduced to a two-dimensional one.

\citeasnoun{Pearson1} used the tetrachoric expansion (\ref{eq_tetra}) for small $|\varrho|$, and quadrature for large $|\varrho|$. \citeasnoun{Nicholson}, building on ideas of \citeasnoun{Sheppard}, worked with a two-parameter function, denoted $V$-function. \citeasnoun{Owen} introduced the $T$-function (\ref{eq_T1}) which is closely related to Nicholson's $V$-function. For many years, quadrature of the $T$-function was the method of choice for evaluation of the bivariate normal distribution. Numerous authors, e.g. \citeasnoun{Borth}, \citeasnoun{Daley}, \citeasnoun{YouMin}, and \citeasnoun{PT}, have been working on improvements, e.g. by dividing the plane into many regions and choosing specific quadrature methods in each region.

\citeasnoun{SowAsh} applied Gauss-Hermite quadrature to (\ref{eq_ab1}) and Simpson's rule to (\ref{eq_1b1}). \citeasnoun{Drezner} used (\ref{eq_uv2line}) and Gauss quadrature. \citeasnoun{Divgi} relied on polar coordinates and an approximation to the univariate Mills' ratio. \citeasnoun{Vasicek} proposed an expansion which is more suitable for large $|\varrho|$ than the tetrachoric expansion (\ref{eq_tetra}).

\citeasnoun{DW} applied Gauss-Legendre quadrature to (\ref{eq_plackett_ind}) for $|\varrho|\leq 0.8$, and to (\ref{eq_plackett_max}) for $|\varrho|>0.8$. Improvements of their method in terms of accuracy and robustness have been provided by \citeasnoun{Genz} and \citeasnoun{West}.

Most implementations today will rely on variants of the approaches of \citeasnoun{Divgi} or of \citeasnoun{DW}. The method of \citeasnoun{Drezner}, although less reliable, is also very common, mainly because it is proposed in \cite{Hull} and other prevalent books.


\section{Bounds and approximations}
\label{sec_bounds}

Nowadays, high-precision numerical evaluation of the bivariate normal copula is usually available and there is not much need for low-precision approximations anymore, in particular if the mathematics is hidden behind strange constants derived from some optimization procedure.

On the other hand, if the mathematics is general and transparent then the resulting approximations are often rather weak. As an example, take approximations to the multivariate Mills' ratios applied to the bivariate case, cf. \cite{LuLi} and the references therein.

Bounds, if they are not too weak, are more interesting than approximations because they can be used, e.g., for checking numerical algorithms. Moreover, derivation of bounds often provides valuable insight into the mathematics behind the function to be bounded.

In the following we concentrate on bounds and approximations explicitly derived for the bivariate case. Throughout this section we will only consider the case $\varrho>0$, $0<u=v<1/2$. Note that by successively applying, if required, (\ref{eq_uv2line}), (\ref{eq_line2diag}), (\ref{eq_gtrafo}) and (\ref{eq_reflect_uv}), we can always reduce $C(u,v;\varrho)$ to a sum of two terms of that form. Any approximation or bound given for the special case can be translated to an approximation or bound for the general case, with at most twice the absolute error. Note also that for many existing approximations and bounds the diagonal $u=v$ may be considered a worst case, cf. \cite{Willink}.

\citeasnoun{MeeOwen} elaborated on the so-called conditional approach proposed by \citeasnoun{Pearson}. If $(X,Y)$ are bivariate standard normally distributed with correlation $\varrho$ then we can write
\begin{equation}
\Phi_2(h,k;\varrho) = \Phi(h)\cdot \PP(Y\leq k\,|\,X\leq h).
\label{eq_condapproach}
\end{equation}
The distribution of $Y$ conditional on $X=h$ is normal but the distribution of $Y$ conditional on $X\leq h$ is not. Nevertheless, it can be approximated by a normal distribution with the same mean and variance. In terms of the bivariate normal copula, the resulting approximation is
\begin{equation}
C(u,u;\varrho) \approx u\cdot\Phi\left(
\frac{u\cdot\Phi^{-1}(u)+\varrho\cdot\varphi\left(\Phi^{-1}(u)\right)}
{\sqrt{u^2-\varrho^2\cdot\varphi\left(\Phi^{-1}(u)\right)\cdot
\left(u\cdot\Phi^{-1}(u)+\varphi\left(\Phi^{-1}(u)\right)\right)}}
\right).
\label{eq_meeowen}
\end{equation}
The approximation works well for $|\varrho|$ not too large. For $|\varrho|$ large there are alternative approximations, e.g. \cite{AK}.

The simpler of the two approximations proposed by \citeasnoun{CoxWermuth} replaces the second factor in (\ref{eq_condapproach}) by the mean of the conditional distribution (the more complicated approximation adds a term of second order). In terms of the bivariate normal copula, the resulting approximation is
\[
C(u,u;\varrho) \approx u\cdot\Phi\left(\sqrt{\frac{1-\varrho}{1+\varrho}}
\cdot\frac{u\cdot\Phi^{-1}(u)+\varrho\cdot\varphi\left(\Phi^{-1}(u)\right)}{(1+\varrho)\cdot u}
\right).
\]

\citeasnoun{Mallows} gave two approximations to Owen's $T$-function (\ref{eq_T1}). In terms of the bivariate normal copula, the simpler one reads
\[
C(u,u;\varrho) \approx 2u\cdot
\Phi\left(\sqrt{\frac{1-\varrho}{1+\varrho}}\cdot\left(\Phi^{-1}\left(\frac{u}{2}+\frac{1}{4}\right)-\Phi^{-1}\left(\frac{3}{4}\right)\right)\right).
\]

Further approximations were derived by \citeasnoun{Cadwell} and \citeasnoun{Polya}.

There are not too many bounds available in the literature. For $\varrho>0$ there are, of course, the trivial bounds (\ref{eq_min}) and (\ref{eq_max}):
\[
u^2 \leq C(u,u;\varrho) \leq u
\]
The upper bound given by \citeasnoun{Polya} is just the one above. His lower bound is too weak (even negative) on the diagonal. A recent overview on known bounds, and derivation of some new ones, is provided by \citeasnoun{Willink}. We will present some of his bounds below in more general context.

\begin{theorem}
Let $\varrho \geq 0$ and $0\leq u \leq 1/2$. Then $C(u,u;\varrho)$ is bounded as follows, where $g(u;\varrho)$ is defined as in (\ref{eq_g}):
\begin{align}
C(u,u;\varrho) & \geq u\cdot g(u;\varrho)\label{eq_bound_lower1}\\
C(u,u;\varrho) & \leq u\cdot g(u;\varrho)\cdot 2\label{eq_bound_upper1}
\end{align}
The lower bound (\ref{eq_bound_lower1}) is tight for $\varrho=0$ or $u=0$. The maximum error of (\ref{eq_bound_lower1}) equals $1/4$ and is obtained for $u=1/2$, $\varrho=1$. The upper bound (\ref{eq_bound_upper1}) is tight for $\varrho=1$ or $u=0$. The maximum error of (\ref{eq_bound_upper1}) equals $1/4$ and is obtained for $u=1/2$, $\varrho=0$.
\label{thm_bound1}
\end{theorem}

A proof of Theorem \ref{thm_bound1} is provided in Section \ref{proof_bound12}, together with a proof of the following refinement:

\begin{theorem}
Let $\varrho > 0$ and $0\leq u \leq 1/2$. Then $C(u,u;\varrho)$ is bounded as follows, where $g(u;\varrho)$ is defined as in (\ref{eq_g}):
\begin{align}
C(u,u;\varrho) & \geq u\cdot g(u;\varrho)\cdot \left(1+\frac{2}{\pi}\arcsin(\varrho)\right)\label{eq_bound_lower2}\\
C(u,u;\varrho) & \leq u\cdot g(u;\varrho)\cdot (1+\varrho)\label{eq_bound_upper2}
\end{align}
These bounds are the optimal ones of the form $u\cdot g(u)\cdot a(\varrho)$. They are tight for $\varrho=0$, $\varrho=1$, or $u=0$. The maximum error of (\ref{eq_bound_upper2}) is obtained for $u=1/2$, $\varrho = \sqrt{1-\frac{4}{\pi^2}} \approx 0.7712$, the value being
\[
\frac{1}{4}\left(\sqrt{1-\frac{4}{\pi^2}}-\frac{2}{\pi}\arcsin\left(\sqrt{1-\frac{4}{\pi^2}}\right)\right) \approx 0.05263.
\]
The lower bound (\ref{eq_bound_lower2}) is tight for $u=1/2$.
\label{thm_bound2}
\end{theorem}

The bounds (\ref{eq_bound_lower1}) and (\ref{eq_bound_upper2}) have been discussed, without explicit computation of the maximum error, by \citeasnoun{Willink}. The maximum error of (\ref{eq_bound_lower2}) is difficult to grasp analytically. Numerically, the error always stays below $0.006$.

An alternative upper bound is given by the following theorem, the proof of which is provided in Section \ref{proof_bound3}:

\begin{theorem}
Let $\varrho > 0$ and $0\leq u \leq 1/2$. Then
\[
C(u,u;\varrho) \leq 2u\cdot g\left(\frac{u}{2};\varrho\right).
\]
The bound is tight for $\varrho=0$, $\varrho=1$, or $u=0$. The maximum error is obtained for $u=1/2$,
\[
-\frac{1+\varrho}{2\pi\cdot\Phi^{-1}(\frac{1}{4})} = \varphi\left(\sqrt{\frac{1-\varrho}{1+\varrho}}\cdot\Phi^{-1}(\textstyle\frac{1}{4})\right),
\]
i.e., $\varrho\approx 0.5961$, the value being approx. $0.015$. 
\label{thm_bound3}
\end{theorem}

It is also possible to derive good approximations to $C(u,u;\varrho)$ by considering the family $C(u,u;\varrho) \approx u\cdot g(u;\varrho) \cdot(a(\varrho)+b(\varrho)u)$. In particular, the choice
\[
a(\varrho) := 1+\varrho, \qquad b(\varrho) := 2\cdot\left(1+\frac{2}{\pi}\arcsin(\varrho)-(1+\varrho)\right) = \frac{4}{\pi}\arcsin(\varrho) - 2\varrho
\]
is attractive because the resulting approximation
\begin{equation}
C(u,u;\varrho) \approx u\cdot g(u;\varrho) \cdot\left(1+\varrho+\left(\frac{4}{\pi}\arcsin(\varrho) - 2\varrho\right)u\right)
\end{equation}
is tight for $\varrho=0$, $\varrho=1$, $u=0$, or $u=1/2$, and for $u\longrightarrow 0^+$ it has the same asymptotic behaviour as (\ref{eq_bound_upper2}). By visual inspection we may conjecture that it is even an upper bound, with an error almost cancelling the error of the lower bound (\ref{eq_bound_lower2}) most of the time. Consequently, an even better approximation (but not a lower bound, for $\varrho$ large) is given by
\begin{equation}
C(u,u;\varrho) \approx u\cdot g(u;\varrho) \cdot\left(1+\frac{\varrho}{2}+\left(\frac{2}{\pi}\arcsin(\varrho) - \varrho\right)u\right).\label{eq_goodapprox}
\end{equation}
Again, (\ref{eq_goodapprox}) is tight for $\varrho=0$, $\varrho=1$, $u=0$, or $u=1/2$. Numerically, the absolute error always stays below 0.0006. Hence (\ref{eq_goodapprox}) is comparable in performance with (\ref{eq_meeowen}), and much better for $\varrho$ large.


\section{Measures of concordance}
\label{sec_concordance}

In the study of dependence between (two) random variables, properties and measures that are scale-invariant, i.e., invariant under strictly increasing transformations of the random variables, can be expressed in terms of the (bivariate) copula of the random variables. Among these are the so-called measures of concordance, in particular Kendall's tau, Spearman's rho, Blomqvist's beta and Gini's gamma. For background and general definitions and properties we refer to Section 5 of \cite{Nelsen}. In this section we will provide formulas for measures of concordance for the bivariate normal copula, depending on the correlation parameter $\varrho$.

Blomqvist's beta follows immediately from (\ref{eq_center}):
\begin{align}
\beta(\varrho) & := 4\cdot C\left(\textstyle\frac{1}{2},\frac{1}{2};\varrho\right) - 1\nonumber\\
& \phantom{:}= \frac{2}{\pi}\arcsin(\varrho)
\label{eq_blomqvist}
\end{align}
For the bivariate normal copula, Kendall's tau equals Blomqvist's beta:
\begin{align}
\tau(\varrho) & := 4\int_0^1 \int_0^1 C(u,v;\varrho) \;d C(u,v;\varrho) - 1\nonumber\\
& \phantom{:}= 1 - 4\int_0^1 \int_0^1
\left( \frac{\partial}{\partial v} C(u,v;\varrho)\right)
\left( \frac{\partial}{\partial u} C(u,v;\varrho)\right) \;du \;dv\nonumber\\
& \phantom{:}= \frac{2}{\pi}\arcsin(\varrho)
\label{eq_kendall}
\end{align}
For a proof of (\ref{eq_blomqvist}) and (\ref{eq_kendall}) cf. Section 5.3.2 of \cite{EFM}. Both Blomqvist's beta and Kendall's tau can be generalized to (copulas of) elliptical distributions, cf. \cite{LMS}. This is not the case for Spearman's rho, cf. \cite{HL}, which is given by:
\begin{align}
\varrho_S(\varrho)
& := 12\int_0^1 \int_0^1 C(u,v;\varrho)-uv\;du\;dv\nonumber\\
& \phantom{:}= 12\int_0^1 \int_0^1 C(u,v;\varrho) \;du\;dv - 3\nonumber\\
& \phantom{:}= \frac{6}{\pi} \arcsin\left(\frac{\varrho}{2}\right)
\label{eq_spearman}
\end{align}
For proofs of (\ref{eq_spearman}) cf. \cite{Kruskal} or Section 5.3.2 of \cite{EFM}.

Gini's gamma for the bivariate normal copula is given as follows:
\begin{align}
\gamma(\varrho) & := 4\left(\int_0^1 C(u,u;\varrho)\;du + \int_0^1 C(u,1-u;\varrho)\;du -\frac{1}{2}\right)\nonumber\\
& \phantom{:}= 4\left(\int_0^1 C(u,u;\varrho)\;du + \int_0^1 u - C(u,u;-\varrho)\;du -\frac{1}{2} \right)\nonumber\\
& \phantom{:}= \frac{2}{\pi}\left(\arcsin\left(\frac{1+\varrho}{2}\right)-\arcsin\left(\frac{1-\varrho}{2}\right)\right)\label{eq_gini1}\\
& \phantom{:} = \frac{4}{\pi} \left(\arcsin\left(\frac{\sqrt{1+\varrho}}{2}\right) - \arcsin\left(\frac{\sqrt{1-\varrho}}{2}\right)\right)\label{eq_gini2}\\
& \phantom{:} = \frac{4}{\pi} \arcsin\left(\frac{1}{4}\left(\sqrt{(1+\varrho)(3+\varrho)}-\sqrt{(1-\varrho)(3-\varrho)}\right)\right)\label{eq_gini3}
\end{align}
I have not been able to find proofs in the literature, hence they will be provided in Section \ref{proof_gini}.

Equation (\ref{eq_gini3}) can be inverted which may be useful for estimation of $\varrho$ from an estimate for $\gamma(\varrho)$:
\begin{equation}
\varrho = \sin\left(\gamma(\varrho)\cdot\frac{\pi}{4}\right)\sqrt{3-\tan\left(\gamma(\varrho)\cdot\frac{\pi}{4}\right)}
\end{equation}

Finally, we propose a new measure of concordance, similar to Gini's gamma but based on the lines $u=1/2$, $v=1/2$ instead of the diagonals. For a bivariate copula $C$ it is defined by:
\begin{align}
\tilde{\gamma}(C(\cdot,\cdot)) & := 4\left(\int_0^1 C\left(u,\textstyle\frac{1}{2}\right)-\frac{u}{2}\;du + \int_0^1 C\left(\textstyle\frac{1}{2},v\right)-\frac{v}{2}\;dv\right)\\
& \phantom{:}= 4\left(\int_0^1 C\left(u,\textstyle\frac{1}{2}\right)\;du + \int_0^1 C\left(\textstyle\frac{1}{2},v\right)\;dv - \frac{1}{2}\right)\nonumber
\end{align}
For the bivariate normal copula we obtain
\begin{equation}
\tilde{\gamma}(\varrho) := \tilde{\gamma}(C(\cdot,\cdot;\varrho))
= \frac{4}{\pi} \arcsin\left(\frac{\varrho}{\sqrt{2}}\right).
\end{equation}
A proof is given implicitly in Section \ref{proof_gini}.


\section{Univariate distributions}

In this section we will discuss two univariate distributions being closely related to the bivariate normal copula (or distribution).

\subsection{The skew-normal distribution}
\label{sec_skew}

A random variable $X$ on $\RR$ is \emph{skew-normally distributed} with skewness parameter $\lambda\in\RR$ if it has a density function of the form
\begin{equation}
f_{\lambda}(x) = 2\varphi(x)\Phi(\lambda x).
\end{equation}
The skew-normal distribution was introduced by \citeasnoun{OHL} and studied and made popular by \citename{Az1985} \citeyear{Az1985,Az1986}. Its cumulative distribution function is given by
\begin{equation}
\PP(X\leq x) = \int_{-\infty}^x 2\varphi(x)\Phi(\lambda x) = 2\int_0^{\Phi(x)} \Phi\left(\lambda\Phi^{-1}(t)\right)\;dt.
\end{equation}
In the light of (\ref{eq_C_diag_int}) and (\ref{eq_diag2line}), cf. also \cite{AzCa}, we find
\begin{align}
\PP(X\leq x) & = 2\Phi_2\left(x,0;-\frac{\lambda}{\sqrt{1+\lambda^2}}\right)\\
& =
\begin{cases}
\Phi_2\left(x,x;\displaystyle\frac{1-\lambda^2}{1+\lambda^2}\right), & \lambda \geq 0,\\
1-\Phi_2\left(-x,-x;\displaystyle\frac{1-\lambda^2}{1+\lambda^2}\right), & \lambda \leq 0.
\end{cases}
\end{align}
In particular, the bounds given in Theorem \ref{thm_bound2} can be applied.

\subsection{The Vasicek distribution}
\label{sec_vasicek}

A random variable $P$ on the interval $[0,1]$ is \emph{Vasicek distributed} with parameters $p\in (0,1)$ and $\varrho\in (0,1)$ if $\Phi^{-1}(P)$ is normally distributed with mean
\begin{equation}
\label{eq_vasicek_mu}
\EE(\Phi^{-1}(P)) = \frac{\Phi^{-1}(p)}{\sqrt{1-\varrho}}
\end{equation}
and variance
\begin{equation}
\label{eq_vasicek_sigma}
\VV(\Phi^{-1}(P)) = \frac{\varrho}{1-\varrho}.
\end{equation}
In Section \ref{sec_model} implicitly it is proved that
\[
\EE(P) = p, \qquad
\EE(P^2) = C(p,p;\varrho),
\]
so that
\[
\VV(P) = \EE(P^2) - \EE(P)^2 = C(p,p;\varrho) - p^2.
\]
Furthermore, we have
\[
\PP(P\leq q) = \PP(\Phi^{-1}(P)\leq\Phi^{-1}(q)) = \Phi\left(\frac{\sqrt{1-\varrho}\cdot\Phi^{-1}(q)-\Phi^{-1}(p)}{\sqrt{\varrho}}\right).
\]
The (one-sided) $\alpha$-Quantile $q_{\alpha}$ of $P$, with $\alpha\in(0,1)$, is therefore given by
\begin{equation}
\label{eq_vasicek_quantile}
q_{\alpha} = \Phi\left(\frac{\sqrt{\varrho}\cdot \Phi^{-1}(\alpha)+\Phi^{-1}(p)}{\sqrt{1-\varrho}}\right).
\end{equation}
In particular, the median of $P$ is simply
\begin{equation}
q_{0.5} = \Phi\left(\frac{\Phi^{-1}(p)}{\sqrt{1-\varrho}}\right) = \Phi\left(\EE(\Phi^{-1}(P))\right).
\end{equation}
The density of $P$ is
\[
\frac{d}{dq} \PP(P\leq q)
= \sqrt{\frac{1-\varrho}{\varrho}}\cdot
\varphi\left(\frac{\sqrt{1-\varrho}\cdot\Phi^{-1}(q)-\Phi^{-1}(p)}{\sqrt{\varrho}}\right)
\cdot\frac{1}{\varphi\left(\Phi^{-1}(q)\right)}.
\]

The distribution is unimodal with the mode at
\[
\Phi\left(\frac{\sqrt{1-\varrho}}{1-2\varrho}\cdot \Phi^{-1}(p)\right)
\]
for $\varrho<0.5$, monotone for $\varrho=0.5$, and U-shaped for $\varrho>0.5$.

Let $\tilde{P}$ be Vasicek distributed with parameters $\tilde{p}$, $\tilde{\varrho}$, and let
\[
\corr\left( \Phi^{-1}(P), \Phi^{-1}(\tilde{P}) \right) = \gamma.
\]
Then
\[
\cov\left( \Phi^{-1}(P), \Phi^{-1}(\tilde{P}) \right)
= \gamma\cdot\sqrt{\frac{\varrho}{1-\varrho}\cdot\frac{\tilde{\varrho}}{1-\tilde{\varrho}}},
\]
\[
\EE\left(P\cdot\tilde{P}\right) = C\left(p,\tilde{p};\gamma\cdot\sqrt{\varrho\cdot \tilde{\varrho}}\right),
\]
and
\[
\cov\left(P,\tilde{P}\right) = C\left(p,\tilde{p};\gamma\cdot\sqrt{\varrho\cdot \tilde{\varrho}}\right) - p\cdot \tilde{p}.
\]

The Vasicek distribution does not offer immediate advantages over other two-parametric continuous distributions on $(0,1)$, such as the beta distribution. Its importance stems from its occurrence as mixing distribution in linear factor models set up as in Section \ref{sec_model}. It is a special case of a probit-normal distribution; it is named after Vasicek who introduced it into credit risk modeling.

For (different) details on the material in this section we refer to \citename{Vasicek_old} \citeyear{Vasicek_old,Vasicek_2002} and \citeasnoun{Tasche}. Estimation of the parameters $p$ and $\varrho$ is also discussed by \citeasnoun{Meyer}.



\begin{appendix}
\section{Proofs}
\label{sec_proofs}

\subsection{Proof of (\ref{eq_abc})}
\label{sec_model}

Let
\begin{align*}
X & = \alpha\cdot Y + \sqrt{1-\alpha^2}\cdot \epsilon,\\
\tilde{X} & = \beta\cdot \tilde{Y} + \sqrt{1-\beta^2}\cdot \tilde{\epsilon},
\end{align*}
where $\alpha,\beta\in(-1,1)\setminus\{0\}$ are parameters and where $Y$, $\tilde{Y}$, $\epsilon$, $\tilde{\epsilon}$ are all standard normal and pairwise independent, except
\[
\gamma := \corr\left(Y,\tilde{Y}\right) = \cov\left(Y,\tilde{Y}\right).
\]
By construction, $X$ and $\tilde{X}$ are standard normal again with
\[
\corr\left(X, \tilde{X}\right) = \cov\left(X, \tilde{X}\right) = \alpha\beta\gamma.
\]
We define indicator variables $Z=Z(X)\in\{0,1\}$, $\tilde{Z}=\tilde{Z}(\tilde{X})\in\{0,1\}$ calibrated to expectation values $u$, $v$:
\begin{align*}
Z = 1 & \qquad:\Longleftrightarrow\qquad X\leq \Phi^{-1}(u),\\
\tilde{Z} = 1 & \qquad:\Longleftrightarrow\qquad \tilde{X}\leq \Phi^{-1}(v)
\end{align*}
Conditional on $(Y=y,\tilde{Y}=\tilde{y})$, $Z$ and $\tilde{Z}$ are independent. We find
\begin{align*}
\PP(Z=1\,|\,Y=y) & = \PP\left(X\leq \Phi^{-1}(u)\,|\,Y=y\right)\\
& = \PP\left(\alpha\cdot y + \sqrt{1-\alpha^2}\cdot\epsilon\leq \Phi^{-1}(u)\right)\\
& = \PP\left(\epsilon\leq\frac{\Phi^{-1}(u)-\alpha\cdot y}{\sqrt{1-\alpha^2}}\right)
= \Phi\left(\frac{\Phi^{-1}(u)-\alpha\cdot y}{\sqrt{1-\alpha^2}}\right).
\end{align*}
Now we define the random variables
\begin{align*}
P & := P(Y) := \Phi\left(\frac{\Phi^{-1}(u)-\alpha\cdot Y}{\sqrt{1-\alpha^2}}\right),\\
\tilde{P} & := \tilde{P}(\tilde{Y}) :=
\Phi\left(\frac{\Phi^{-1}(v)-\beta\cdot Y}{\sqrt{1-\beta^2}}\right).
\end{align*}
We find
\[
u = \EE(Z) = \EE(\EE(Z\,|\,Y)) = \EE(\PP(Z=1\,|\,Y)) = \EE(P),
\qquad v = \EE(\tilde{P})
\]
and
\begin{align*}
\PP\left(Z=1,\tilde{Z}=1\right) & = \PP\left(X\leq \Phi^{-1}(u), \tilde{X}\leq \Phi^{-1}(v)\right)\\
& = \Phi_2\left(\Phi^{-1}(u),\Phi^{-1}(v),\cov(X,\tilde{X})\right)\\
& = \Phi_2\left(\Phi^{-1}(u),\Phi^{-1}(v),\alpha\beta\gamma\right).
\end{align*}
On the other hand,
\begin{align*}
& \PP\left(Z=1,\tilde{Z}=1\right) = \PP\left(Z \tilde{Z} = 1\right) = \EE\left(Z \tilde{Z}\right)\\
& \qquad = \EE\left(\EE(Z \tilde{Z}\,|\,Y,\tilde{Y})\right) = \EE\left(\EE(Z\,|\,Y,\tilde{Y}) \cdot\EE(\tilde{Z}\,|\,Y,\tilde{Y})\right)\\
& \qquad = \EE\left(\EE(Z\,|\,Y) \cdot\EE(\tilde{Z}\,|\,\tilde{Y})\right) = \EE\left(P\cdot\tilde{P}\right)\\
& \qquad = \int_{-\infty}^{\infty} \int_{-\infty}^{\infty} P(x)\cdot \tilde{P}(y)\cdot \varphi_2(x,y,\gamma) \;dx \;dy\\
& \qquad = \int_{-\infty}^{\infty} \int_{-\infty}^{\infty} \Phi\left(\frac{\Phi^{-1}(u)-\alpha\cdot x}{\sqrt{1-\alpha^2}}\right)
\Phi\left(\frac{\Phi^{-1}(v)-\beta\cdot y}{\sqrt{1-\beta^2}}\right)
\varphi_2(x,y,\gamma) \;dx \;dy.
\end{align*}

\subsection{Proof of Theorems \ref{thm_bound1} and \ref{thm_bound2}}
\label{proof_bound12}

We will assume $\varrho$ as fixed and write
\[
C(u) := C(u,u;\varrho), \qquad g(u) := g(u;\varrho).
\]

The upper bound (\ref{eq_bound_upper1}) follows from (\ref{eq_gtrafo}).
Regarding the lower bound (\ref{eq_bound_lower1}) we note that
\begin{align*}
g'(u) &
= \sqrt{\frac{1-\varrho}{1+\varrho}}\cdot\exp\left(\frac{\varrho}{1+\varrho}\cdot\Phi^{-1}(u)^2\right) > 0,\\
g''(u) &
= g'(u)\cdot\frac{2\varrho}{1+\varrho}\cdot\frac{\Phi^{-1}(u)}{\varphi(\Phi^{-1}(u))} < 0.
\end{align*}
Hence $g$ is increasing and concave on $(0,1/2)$ and we conclude
\[
C(u) = 2\int_0^u g(t)\;dt \geq 2\cdot\frac{1}{2}\cdot u\cdot g(u) = u\cdot g(u).
\]
Now we define
\[
D_a(u) := u\cdot g(u)\cdot a
\]
with $a\in[1,2]$, $u\in [0,1/2]$. We start by noting that
\[
C'(u) = 2g(u) > 0, \qquad
C''(u) = 2g'(u) > 0,
\]
hence $C$ is increasing and convex on $(0,1/2)$. Furthermore,
\begin{align*}
D_a'(u) & = a\cdot\left(g(u)+u\cdot g'(u)\right) > 0,\\
D_a''(u) & =
a\cdot\left(2\cdot g'(u)+u\cdot g''(u)\right)\\
& = 2a\cdot g'(u)\cdot\left(1+\frac{1}{1+\varrho}\cdot\frac{u\cdot\Phi^{-1}(u)}{\varphi\left(\Phi^{-1}(u)\right)}\right)\\
& = 2a\cdot g'(u)\cdot\left(1-\frac{1}{1+\varrho}\cdot H\left(-\Phi^{-1}(u)\right)\right)
\end{align*}
with $H(x)=x\cdot R(x)$, where
\[
R(x) = \frac{1-\Phi(x)}{\varphi(x)}
\]
is Mills' ratio. \citeasnoun{Pinelis} has shown that $H'(x)>0$ for $x>0$, $H(0)=0$, and $\lim_{x\longrightarrow\infty}H(x)=1$. Hence $D_a(u)$ is increasing and convex on $(0,1/2)$ as well.

For $u\in(0,1/2)$, $C'(u) = D_a'(u)$ is equivalent with
\[
f(u) := \frac{u\cdot g'(u)}{g(u)} = \frac{2-a}{a}, \qquad \text{ or }\qquad a = \frac{2}{1+f(u)}.
\]
We will show that $f$ is strictly increasing on $(0,1/2)$. We have
\[
f(u) = \lambda\cdot\frac{\Phi\left(\Phi^{-1}(u)\right)\cdot\varphi\left(\lambda\cdot\Phi^{-1}(u)\right)}{\varphi\left(\Phi^{-1}(u)\right)\cdot\Phi\left(\lambda\cdot\Phi^{-1}(u)\right)}
= \lambda\cdot F_{\lambda}\left(-\Phi^{-1}(u)\right)
\]
with $\lambda = \sqrt{\frac{1-\varrho}{1+\varrho}}\in [0,1]$ and
\[
F_{\lambda}(x) := \frac{R(x)}{R(\lambda\cdot x)}, \qquad x\geq 0.
\]
We find
\begin{align*}
F_{\lambda}'(x) & = \frac{R'(x)\cdot R(\lambda\cdot x)-R(x)\cdot R'(\lambda\cdot x)}{R(\lambda\cdot x)^2}\\
& = F_{\lambda}(x)\cdot\left(\frac{R'(x)}{R(x)}-\lambda\cdot\frac{R'(\lambda\cdot x)}{R(\lambda\cdot x)}\right) < 0
\end{align*}
for $\lambda < 1$. Here we have used that $F_{\lambda}(x) > 0$ and that the function
\[
y \mapsto y\cdot\frac{R'(y)}{R(y)}
\]
is strictly decreasing on $(0,\infty)$, cf. \cite{Baricz}. We conclude
\[
f'(u) = -\lambda\cdot \frac{F_{\lambda}'\left(-\Phi^{-1}(u)\right)}{\varphi\left(-\Phi^{-1}(u)\right)} > 0.
\]
Furthermore, we find
\[
f\left(\textstyle\frac{1}{2}\right)=\lambda\cdot\frac{R(0)}{R(0)} = \sqrt{\frac{1-\varrho}{1+\varrho}}
\]
and
\begin{align*}
\lim_{u\longrightarrow 0^+} f(u) & = \lim_{u\longrightarrow 0^+} \frac{g'(u)+u\cdot g''(u)}{g'(u)}\\
& = \lim_{u\longrightarrow 0^+} 1 + \frac{2\varrho}{1+\varrho}\cdot\frac{u\cdot\Phi^{-1}(u)}{\varphi(\Phi^{-1}(u))}\\
& = \lim_{u\longrightarrow 0^+} 1 - \frac{2\varrho}{1+\varrho}\cdot H\left(-\Phi^{-1}(u)\right)\\
& = 1 - \frac{2\varrho}{1+\varrho} = \frac{1-\varrho}{1+\varrho}.
\end{align*}

We have $C(0)=D_a(0)=0$, $C'(0)=D_a'(0)=0$, and
\[
\lim_{u\longrightarrow 0^+} \frac{D_a(u)}{C(u)} =
\lim_{u\longrightarrow 0^+} \frac{D_a'(u)}{C'(u)} =
\lim_{u\longrightarrow 0^+} \frac{a}{2}\left(1+f(u)\right) = \frac{a}{1+\varrho}.
\]
By standard calculus we conclude that
\begin{itemize}
\item For $a\geq 1+\varrho$ we have $D_a'(u)\geq C'(u)$, and hence $D_a(u)\geq C(u)$, for all $u\in[0,1/2]$;
\item For $a\leq 2\cdot\left(1+\sqrt{\frac{1-\varrho}{1+\varrho}}\right)^{-1}$ we have $D_a'(u)\leq C'(u)$, and hence $D_a(u)\leq C(u)$, for all $u\in[0,1/2]$;
\item For
\[
a\in\left(2\cdot\left(1+\sqrt{\frac{1-\varrho}{1+\varrho}}\right)^{-1}, 1+\varrho\right)
\]
there exists $u_0\in (0,1/2)$ with $D_a'(u)<C'(u)$ for $u\in (0,u_0)$, and $D_a'(u)>C'(u)$ for $u\in (u_0,1/2)$. Consequently, the best lower bound for $C$ of the form $D_a$ is obtained if $C(1/2)=D_a(1/2)$, i.e., $a=1+\frac{2}{\pi}\arcsin(\varrho)$. Moreover, the upper bound $D_a$ with $a=1+\varrho$ can not be improved.
\end{itemize}

The maximum error of $D_a$ with $a=1+\varrho$ is attained if
\[
\frac{d}{d\varrho}\left[ D_a(1/2) - C(1/2)\right] = \frac{1}{4}-\frac{1}{2\pi\sqrt{1-\varrho^2}} = 0,
\]
which is equivalent with $\varrho=\sqrt{1-\frac{4}{\pi^2}}$.

\subsection{Proof of Theorem \ref{thm_bound3}}
\label{proof_bound3}

We will assume $\varrho$ as fixed and write
\[
C(u) := C(u,u;\varrho), \qquad g(u) := g(u;\varrho).
\]
We have
\[
C(u) = 2\cdot \int_0^u g(t)\;dt = 2u\cdot g(v(u))
\]
with
\[
v(u):=v(u;\varrho)\leq u.
\]
Since $g$ is concave and increasing, we even know that
\[
v(u) \leq \frac{u}{2}
\]
and hence
\[
C(u) = 2u\cdot g(v(u)) \leq 2u\cdot g\left(\frac{u}{2}\right).
\]
Moreover, for the same reason we have
\begin{align*}
\frac{d}{du} \left(2u\cdot g\left(\frac{u}{2}\right) - C(u)\right)
& = 2\left(\frac{u}{2}\cdot g'\left(\frac{u}{2}\right)-\left(g(u)-g\left(\frac{u}{2}\right)\right)\right) \geq 0,
\end{align*}
and hence, for $\varrho$ fixed, the maximum error is obtained for $u=1/2$, the value being
\[
\Phi\left(\sqrt{\frac{1-\varrho}{1+\varrho}}\cdot\Phi^{-1}\left(\frac{1}{4}\right)\right)-\frac{1}{4}-\frac{1}{2\pi}\arcsin(\varrho).
\]
Derivation of the above expression with respect to $\varrho$ gives the result.

Note that by (\ref{eq_g_inv}) we can write
\[
v(u;\varrho) = g\left(\frac{C(u)}{2u};-\varrho\right)
\]
Unfortunately, for large $\varrho$, the function $v$ is not convex, and the approximation
\[
2u\cdot g\left(\frac{u}{2}\right)\cdot\frac{C(\frac{1}{2})}{g(\frac{1}{4})}
\]
is not an upper bound for $C(u)$.


\subsection{Proof of (\ref{eq_gini1}), (\ref{eq_gini2}), (\ref{eq_gini3})}
\label{proof_gini}
In a first step, using (\ref{eq_plackett_ind}) we find:
\begin{align*}
\int_0^1 C(u,u;\varrho) \;du
& = \int_0^1 u^2 + \frac{1}{2\pi} \int_0^{\varrho} \frac{1}{\sqrt{1-r^2}}\exp\left(-\frac{\Phi^{-1}(u)^2}{1+r}\right)\;dr\;du\\
& = \frac{1}{3} + \frac{1}{2\pi} \int_0^{\varrho} \frac{1}{\sqrt{1-r^2}} \int_{-\infty}^{\infty} \varphi(v) \exp\left(-\frac{v^2}{1+r}\right) \;dv\;dr\\
& = \frac{1}{3} + \frac{1}{2\pi} \int_0^{\varrho} \frac{1}{\sqrt{1-r^2}} \int_{-\infty}^{\infty} \varphi(s)\sqrt{\frac{1+r}{3+r}}\;ds\;dr\\
& = \frac{1}{3} + \frac{1}{2\pi} \int_0^{\varrho} \frac{1}{\sqrt{(1-r)(3+r)}}\;dr\\
& = \frac{1}{3} + \frac{1}{2\pi} \int_{\frac{1}{2}}^{\frac{1+\varrho}{2}}\frac{1}{\sqrt{1-r^2}}\;dr\\
& = \frac{1}{3} + \frac{1}{2\pi} \left(\arcsin\left(\frac{1+\varrho}{2}\right)-\arcsin\left(\frac{1}{2}\right)\right)\\
& = \frac{1}{4} + \frac{1}{2\pi} \arcsin\left(\frac{1+\varrho}{2}\right).
\end{align*}
We conclude that
\begin{align*}
\gamma(\varrho) & = 4\left(\int_0^1 C(u,u;\varrho)\;du + \int_0^1 u - C(u,u;-\varrho)\;du -\frac{1}{2} \right)\nonumber\\
& = 4\left(\frac{1}{4} + \frac{1}{2\pi} \arcsin\left(\frac{1+\varrho}{2}\right)+\frac{1}{2}-\frac{1}{4} - \frac{1}{2\pi} \arcsin\left(\frac{1-\varrho}{2}\right)-\frac{1}{2}\right)\nonumber\\
& = \frac{2}{\pi}\left(\arcsin\left(\frac{1+\varrho}{2}\right)-\arcsin\left(\frac{1-\varrho}{2}\right)\right).
\end{align*}
In a similar way, using again (\ref{eq_plackett_ind}), we can compute
\begin{align*}
\int_0^1 C\left(u,\textstyle\frac{1}{2};\varrho\right)\;du
& = \int_0^1 \frac{u}{2} + \frac{1}{2\pi}\int_0^{\varrho} \frac{1}{\sqrt{1-r^2}}\exp\left(-\frac{\Phi^{-1}(u)^2}{2(1-r^2)}\right)\;dr\;du\\
& = \frac{1}{4} + \frac{1}{2\pi} \arcsin\left(\frac{\varrho}{\sqrt{2}}\right),
\end{align*}
which, using (\ref{eq_diag2line}), leads to
\[
\int_0^1 C(u,u;\varrho)\;du = 2\cdot \int_0^1 C\left(u,\textstyle\frac{1}{2};-\sqrt{\frac{1-\varrho}{2}}\right)\;du = \frac{1}{2} - \frac{1}{\pi} \arcsin\left(\frac{\sqrt{1-\varrho}}{2}\right).
\]
We obtain alternative formulas for Gini's gamma, the second one using the addition theorem for the $\arcsin$ function:
\begin{align*}
\gamma(\varrho) & = \frac{4}{\pi} \left(\arcsin\left(\frac{\sqrt{1+\varrho}}{2}\right) - \arcsin\left(\frac{\sqrt{1-\varrho}}{2}\right)\right)\\
& = \frac{4}{\pi} \arcsin\left(\frac{1}{4}\left(\sqrt{(1+\varrho)(3+\varrho)}-\sqrt{(1-\varrho)(3-\varrho)}\right)\right)
\end{align*}

\end{appendix}

\end{document}